\documentclass[11pt]{article}
\usepackage{amssymb,amsmath,euscript}
\date{June 2, 2009}
\newtheorem{proposition}{Proposition}[section]

\newtheorem{lemma}[proposition]{Lemma}
\newtheorem{theorem}[proposition]{Theorem}
\newtheorem{definition}[proposition]{Definition}
\newtheorem{corollary}[proposition]{Corollary}

\def\l{{\langle}}
\def\r{\rangle}

\def\R{{\mathbb R}}

\def\a{\alpha}

\def\De{\Delta}

\def\ga{\gamma}
\def\Ga{\Gamma}
\def\ep{\varepsilon}
\def\eps{\varepsilon}

\def\Re {{\rm Re}\,}

\def\E{{\mathbb E}}
\def\P{{\mathbb P}}

\makeatletter \@addtoreset{equation}{section} \makeatother
\newcommand {\qed}%
{%
    {}\hfill
    {}\hfill
    {$\square $}%
    \vspace {0.3cm}%
    \pagebreak [2]%
    \par
}%
\newenvironment{proof}[1]{%
    \vspace{0.3cm}%
    \pagebreak [2]%
    \par%
    \noindent {\bf  Proof~#1\ }}{\qed}%
\newenvironment{example}{%
    \vspace{0.3cm} \pagebreak [2]%
    \par%
    \refstepcounter{proposition}%
    \noindent%
    {\bf  Example~\theproposition\ }}{\ }%
\newenvironment{remark}{%
    \vspace{0.3cm} \pagebreak [2]%
    \par%
    \refstepcounter{proposition}
    \noindent%
   {\bf Remark~\theproposition\  }}{\ }%
\begin{document}
%


\title {Uniform Modulus of Continuity of Random Fields}
\author{Yimin Xiao
\thanks{Research supported in part by the NSF grant DMS-0706728.}
\\ Michigan State University}

\maketitle

\begin{abstract}
A sufficient condition for the uniform modulus of continuity of a random
field $X = \{X(t), t \in \R^N\}$ is provided. The result is applicable to
random fields with heavy-tailed distribution such as stable random fields.
\end{abstract}

{\sc Running head}:  Uniform Modulus of Continuity of Random Fields\\

{\sc 2000 AMS Classification numbers}: 60G60, 60G52, 60G17, 60G18.\\


\section{Introduction}

Sample path continuity and H\"older regularity of stochastic processes
have been studied by many authors. The celebrated theorems of
Kolmorogov [cf. Khoshnevisan (2002)] and Garsia (1972) provide uniform
modulus of continuity for rather general stochastic processes
and random fields. For Gaussian processes, a powerful chaining argument
leads to much deeper results; see the recent books of Talagrand (2006),
Marcus and Rosen (2006) and Adler and Taylor (2007).
Some of the results on uniform modulus of continuity have been extended
to non-Gaussian processes provided the tail
probabilities of the increments of the process decay fast enough [see
Cs\'{a}ki and Cs\"{o}rg\H{o} (1992), Kwapie\'n and Rosi\'nski (2004)].

In recent years, there has been increasing interest in sample path
continuity of stochastic processes with heavy-tailed distributions such as
certain infinitely divisible processes including particularly stable processes.
See Samorodnitsky and Taqqu (1994),
Marcus and Rosi\'nski (2005) and the references therein for more
information. Since many of such processes do not have finite second moment,
the aforementioned results on the uniform modulus of continuity do not apply.
There have been only a few results on uniform modulus of continuity for
stable processes; see Bernard (1970), Takashima (1989), K\^ono and Maejima
(1991a, 1991b), Bierm\'e and Lacaux (2009), Ayache, Roueff and Xiao
(2007, 2009).

The objective of this paper is to study uniform modulus of continuity of a
real-valued random field $X=\{X(t), t \in \R^N\}$ with heavy-tailed distributions.
In Section 2, we
modify the chaining argument to prove a general result on uniform modulus
of continuity for $X$. Compared with Kolmogorov's continuity theorem,
Dudley's entropy theorem or the theorems in Cs\'{a}ki and Cs\"{o}rg\H{o}
(1992), our result does not assume exponential-type
tail probabilities nor higher moments, hence it can be applied to wider
classes of random fields. In Section 3, we introduce the notion of maximal
moment index for a sequence of random variables and study its basic properties.
We show that the maximal moment index is closely related to the uniform
H\"older exponent of the random field $X$. In Section 4, we apply the
results in Sections 2 and 3 to stable random fields.

Throughout this paper,  we will use $K$ to denote an unspecified positive
constant which may differ in each occurrence. Some specific constants
will be denoted by $K_{1}, K_2, \ldots$.

\bigskip

\noindent{\bf Acknowledgment}\ The author thanks Professors Michael B. Marcus
and Jan Rosi\'nski for stimulating discussions. He is grateful to the anonymous
referee for his/her comments and suggestions which have led to significant
improvement of the manuscript.

This paper was finished during the authors visit to the Statistical \&
Applied Mathematical Sciences Institute (SAMSI).
He thanks the staff of SAMSI for their support and the good
working conditions.

\section{A general result on uniform modulus of continuity}
\label{sec:MC}

Let $X = \{X(t), t \in T\}$ be a real-valued random field on a compact
metric space $(T, d)$. In this section we establish a general result on
uniform modulus of continuity of $X$. Similar to the proofs of Kolmogorov's
continuity theorem and Dudley's entropy theorem [see Khoshnevisan (2002),
Talagrand (2006), Marcus and Rosen (2006) or Adler and Taylor (2007)],
our method is also based
on the powerful chaining argument. However, in order to be able to apply
the method to random fields with heavy-tailed distributions such as stable
random fields, we will make two modifications. One is that we will rely
more on the geometry of the parameter set $T$ and choose approximating
chains more carefully so that, at the $n$-th step of the chain, there
are $\kappa_0\,N_n$ pairs of random variables, where $N_n$ is the
$2^{-n}$-entropy number of $T$ and $\kappa_0$ is a constant, instead
of $N_nN_{n-1}$ pairs as in Talagrand (2006) or Adler and Taylor (2007).
This is important for deriving correct rate function for the modulus of
continuity for processes with heavy-tailed distributions.
The second modification is, instead of trying to control the tail probabilities of
the increment of $X(t)$ over every step of the chain as in the
aforementioned references, we control the $\ga$-th moments of the
maximum increments over the steps. Here $\gamma \in (0, 1]$ is a constant.
This latter argument extends an idea of K\^ono and Maejima (1991b)
who studied the modulus of continuity of an $H$-self-similar $\alpha$-stable
process with stationary increments and satisfying $\alpha \in (1, 2)$
and $\frac 1 \alpha < H < 1$.

We assume that there is a sequence
$\{D_n, n \ge 1\}$ of finite subsets of $T$ satisfying the following conditions:
\begin{itemize}
\item[(i).] There exists a positive integer $\kappa_0$ depending only on $(T, d)$
such that for every $\tau_{n+1} \in D_{n+1}$, the number of
points $\tau_n \in  D_{n}$ satisfying $d(\tau_{n+1}, \tau_n) \le
2^{-n}$ is at most $\kappa_0$. Such a point $\tau_n$ is called a $D_{n}$-neighbor
of $\tau_{n+1}$. We denote by $O_n(\tau_{n+1})$ the set of all $D_{n}$-neighbors
of $\tau_{n+1}$.
\item[(ii).] [Chaining property] For every $s, t \in T$ with $d(s,t) \le 2^{-n}$,
there exist two sequences
$\{\tau_{p}(s), p \ge n\}$ and $\{\tau_{p}(t), p \ge n\}$ such that $\tau_{n}(s)
= \tau_{n}(t)$ (i.e., the two chains have the same starting point) and, for every
$p \ge n$, both $\tau_{p}(s)$ and $\tau_{p}(t)$ are in $D_p$,
$$
d(\tau_{p}(t), \, t) \le  2^{-p}, \ \ \  \hbox{  } \ \ \ \ d(\tau_{p}(s),\,  s)
\le  2^{-p}
$$
and
$$\tau_{p}(t)\in O_p\big(\tau_{p+1}(t)\big),\ \ \ \ \hbox{  \ }\ \ \
\tau_{p}(s)\in O_p\big(\tau_{p+1}(s)\big).
$$
If, in addition,  $s \in D := \bigcup_{k=1}^\infty D_k$ (or $t \in D$), then there
exists an integer $q \ge n$ such that
$\tau_{p}(s) = s$  (or $\tau_{p}(t) = t$) for all $p \ge q$.
\end{itemize}

Condition (ii) implies that, for every $n$, $T$ is covered by $d$-balls
with radius $2^{-n}$ and centers in $D_n$; and the collection of all
the points in $\{D_n, n \ge 1\}$ is dense in $T$, i.e.,
$\overline{\bigcup_{n=1}^\infty D_n} = T$.

When $T$ is a closed interval in $\R^N$ and $d$ is a metric which is
equivalent to the $\ell^\infty$ metric $|s-t|_\infty = \max_{1
\le j \le N} |s_j- t_j|$, it is convenient to choose a sequence
$\{D_n,  n \ge 1\}$ satisfying the above conditions. For example,
if $T= [0, 1]^N$ and $d$ is the $\ell^\infty$ metric, then, for any
integer $n \ge 1$, one can take $D_n$ to be the collection of all
vertices of dyadic cubes of order $n$ which are contained in $T$.
Note that in this case, $\{D_n, n \ge 1\}$ is increasing (i.e.,
$D_n \subset D_{n+1}$ for every $n \ge 1$) and for every $\tau_n
\in D_n$, the number of $D_{n-1}$-neighbors of
$\tau_n$  is at most $\kappa_0:=3^N-1$.


The following is the main result of this section.
\begin{theorem}\label{Thm:MC1}
Let $X = \{X(t), t \in T\}$ be a real-valued random field on a compact
metric space $(T, d)$ and let $\{D_n, n \ge 1\}$ be a sequence of finite subsets of
$T$ satisfying Conditions (i) and (ii).
Suppose $\sigma: \R_+ \to \R_+$ is a nondecreasing continuous function such that
$\sigma(2h) \le K_1 \sigma(h)$  for some constant $K_1 > 0$. If there exist
constants $K_2 > 0$, $\gamma \in (0, 1]$ and $\ep_0 > 0$ such that
$\lim\limits_{h \to 0} \sigma(h)
\big(\log (1/h) \big)^{(1 + \ep_0)/\ga} = 0$ and
\begin{equation}\label{Eq:MomCon1B}
\sum_{p=n}^\infty \E\bigg(\max_{\tau_p \in D_p} \max_{\tau'_{p-1} \in O_{p-1}(\tau_{p})}
\big|X(\tau_p) - X(\tau_{p-1}')\big|^\gamma\bigg)
\le K_2\, \sigma(2^{- n})^\gamma.
\end{equation}
Then $X$ has a continuous version, still denoted by $X$, such that for all $\ep > 0$
\begin{equation}\label{Eq:UM1c}
\lim_{h\to 0+} \,\frac{\sup_{t \in T}\sup_{d(s,t)\le h}
\big|X(t) - X(s)\big|} {\sigma(h)\, \big(\log (1/h)\big)^{(1 + \ep)/\ga }} = 0,
\qquad \hbox{ a.s.}
\end{equation}
\end{theorem}

\begin{proof}\ Given any $s, t \in D = \bigcup_{k=1}^\infty D_k$ with
$d(s,t) \le 2^{-n}$, let $\{\tau_{p}(s), n \le p \le q\}$ and $\{\tau_{p}(t),
n \le p \le q\}$ be the two approximating chains to $s$ and $t$ given
by Condition (ii). The triangle inequality and the fact $\tau_n(s) = \tau_n(t)$
imply
\begin{equation}\label{Eq:Chain1}
\begin{split}
\big|X(t) - X(s)\big| &\le \sum_{p = n+1}^q \big|X(\tau_p(t)) - X(\tau_{p-1}(t))\big|
+ \sum_{p = n+1}^q \big|X(\tau_p(s)) - X(\tau_{p-1}(s))\big|\\
&\le 2 \sum_{p= n+1}^\infty \max_{\tau_p \in D_p} \max_{\tau'_{p-1} \in O_{p-1}(\tau_n)}
\big|X(\tau_p) - X(\tau_{p-1}')\big|.
\end{split}
\end{equation}
It is helpful to note that, for each $p \ge n+1$, the maximum in (\ref{Eq:Chain1})
is taken over at most $\kappa_0\,N_p $ increments, where $N_p$ denotes the cardinality
of $D_p$.

For any integer $n \ge 1$, let
\[
\Delta_n = \sup_{s,t \in D} \sup_{d(s, t) \le 2^{-n}}\big|X(t) - X(s)\big|.
\]
It follows from (\ref{Eq:Chain1}), the elementary inequality
$|a+b|^\ga \le |a|^\ga + |b|^\ga$ ($\ga \in (0, 1]$ and $a, b \in \R$) and (\ref{Eq:MomCon1B})
that
\[
\E\big(\Delta_n^\ga\big) \le 2^\ga \sum_{p= n+1}^\infty
\E\bigg(\max_{\tau_p \in D_p} \max_{\tau'_{p-1} \in O_{p-1}(\tau_p)}
\big|X(\tau_p) - X(\tau_{p-1}')\big|^\ga\bigg)\le 2^\ga\,K_2 \sigma( 2^{-n})^\ga.
\]
Hence Markov's inequality gives that for all integers $n \ge 1$ and real
numbers $u > 0$
\begin{equation}\label{Eq:Chain3}
\P\Big( \Delta_n \ge \sigma( 2^{-n})\, u \Big) \le 2^\ga\,K_2\, u^{-\ga}.
\end{equation}
For any $\ep \in (0, \ep_0)$ and $\eta>0$, by taking
$u = \eta (\log 2^n)^{(1+\ep)/\ga}$, we derive from (\ref{Eq:Chain3}) that
\[
\P\Big( \Delta_n \ge \eta\, \sigma(2^{-n})\, (\log 2^n)^{(1+\ep)/\ga}
\Big) \le \frac{K_3} {n^{1 + \ep}},
\]
which is summable. The Borel-Cantelli lemma implies almost surely
\begin{equation} \label{Eq:Chain5a}
\Delta_n \le \eta\, \sigma(2^{- n})\, (\log 2^n)^{(1 + \ep)/\ga}
\end{equation}
for all $n$ large enough. Therefore, $X(t)$ is a.s. uniformly continuous
in $D = \bigcup_{k=1}^\infty D_k$. Since $D$ is dense in $T$, it is standard
to derive from (\ref{Eq:Chain5a}) that $X$ has a continuous version
(still denoted by $X$) such that for all $\ep > 0$,
\begin{equation} \label{Eq:UM1b}
\lim_{n \to \infty} \frac{\sup_{t \in T} \sup_{d(s, t) \le 2^{-n}}\big|X(t)-X(s)\big|}
{\sigma(2^{- n})\, (\log 2^n)^{(1 + \ep)/\ga}} = 0,
\qquad \hbox{ a.s.}
\end{equation}
By (\ref{Eq:UM1b}), the properties of $\sigma$ and a monotonicity argument
we derive (\ref{Eq:UM1c}). This finishes the proof of Theorem \ref{Thm:MC1}.
\end{proof}

From now on, we will not distinguish $X$ from its continuous version.  The
following result follows directly from Theorem \ref{Thm:MC1} and is often
more convenient to use.

\begin{corollary}\label{Thm:MC1b}
Let $X = \{X(t), t \in T\}$ be a real-valued random field on a compact
metric space $(T, d)$ and let $\{D_n, n \ge 1\}$ be as in Theorem \ref{Thm:MC1}.
Assume that there exist constants $\gamma \in (0, 1]$, $\delta > 0$ and $K>0$  such
that
\begin{equation}\label{Eq:MomCon1}
\E\bigg(\max_{\tau_n \in D_n} \max_{\tau'_{n-1} \in O_{n-1}(\tau_{n}) }
\big|X(\tau_n) - X(\tau_{n-1}')\big|^\gamma\bigg)
\le K\, 2^{-\delta \gamma n}
\end{equation}
for all integers $n \ge 1$. Then for all $\ep > 0$,
\begin{equation}\label{Eq:UM1}
\lim_{h\to 0+} \frac{\sup_{t \in T}\sup_{ d(s,t)\le h} \big|X(t) - X(s)\big|}
{h^\delta\, \big(\log (1/h)\big)^{(1 + \ep)/\ga }} = 0,
\qquad \hbox{ a.s.}
\end{equation}
\end{corollary}

Next we apply Theorem \ref{Thm:MC1} to a random field $X =\{X(t), t \in [0, 1]^N\}$
which may be anisotropic. Corollary \ref{Thm:MC2} below is applicable to anisotropic
stable random fields including linear fractional stable sheets considered in
Ayache, Roueff and Xiao (2007, 2009),
harmonizable fractional stable sheets [cf, Xiao (2006, 2008)], operator-scaling
stable fields with stationary increments in Bierm\'e and Lacaux (2009) and
others infinitely divisible fields.

Given a constant vector $(H_1, \ldots, H_N)\in (0, 1]^N$,
we consider the metric $\rho$ on $\R^N$ defined by:
\begin{equation}\label{Def:rho}
\rho(s, t) = \max_{1 \le j \le N} |s_j - t_j|^{H_j},\qquad \forall\, s,\, t \in \R^N.
\end{equation}
For every $n \ge 1$, define
\begin{equation}\label{Eq:Dn2}
\widetilde D_n = \bigg\{\Big(\frac{k_1} {2^{n/H_1}}, \ldots, \frac{k_N} {2^{n/H_N}}\Big),
\ 1 \le k_j \le \lfloor 2^{n/H_j}\rfloor, \ \forall 1 \le j \le N\bigg\}.
\end{equation}
Then
$$\# \widetilde D_n \le 2^{Qn} \ \ \ \hbox{ and }\ \ \
\overline{\bigcup_{n=1}^\infty \widetilde D_n} = [0, 1]^N.
$$
In the above and in the sequel, $\# D$ denotes the cardinality of
$D$ and $Q = \sum_{j=1}^N \frac 1 {H_j}$. It is easy to see that the
sequence $\{\widetilde D_n, n \ge 1\}$ satisfies Conditions (i) and (ii) on the
compact metric space $([0, 1]^N, \rho)$.

\begin{corollary}\label{Thm:MC2}
Let $X = \{X(t), t \in [0, 1]^N\}$ be a real-valued random field and assume
$\sigma: \R_+\to \R_+$ satisfies the conditions of Theorem \ref{Thm:MC1}. Let
$\{\widetilde{D}_n, n \ge 1\}$ be the sequence defined above. If there exist
constants $\gamma \in (0, 1]$ and $K > 0$ such that
\[
\sum_{p=n}^\infty \E\bigg(\max_{\tau_p \in \widetilde{D}_p}
\max_{\tau'_{p-1}\in O_{p-1}(\tau_p)}
\big|X(\tau_p) - X(\tau_{p-1}')\big|^\gamma\bigg)
\le K\, \sigma\big(2^{-n}\big)^{\ga}
\]
for all integers $n \ge 1$, then for any $\ep > 0$,
\begin{equation}\label{Eq:UM1d}
\lim_{h\to 0+} \sup_{t \in [0,1]^N}\sup_{ |s-t|\le h} \frac{\big|X(t) - X(s)\big|}
{\sigma\big(\sum_{j=1}^N |t_j - s_j|^{H_j}\big)\,\big|
\log\big(\sum_{j=1}^N |t_j - s_j|^{H_j} \big)\big|^{(1+\eps)/\ga}} = 0,
\quad \hbox{a.s.}
\end{equation}
\end{corollary}

\begin{proof}\ Note that the conditions of Theorem  \ref{Thm:MC1} are satisfied.
By (\ref{Eq:UM1c}) we have that for all $\ep > 0$
\begin{equation}\label{Eq:UM3a}
\lim_{h\to 0+} \,\frac{\sup_{t \in [0, 1]^N, \, \rho(s,t)\le h}
\big|X(t) - X(s)\big|} {\sigma(h)\, \big(\log (1/h)\big)^{\frac{1 + \ep} \ga}} = 0,
\qquad \hbox{ a.s.}
\end{equation}
For any $s, t \in [0, 1]^N$, there is an integer $n$ such that $ 2^{-n-1}\le \rho(s, t)
< 2^{-n}$. Then (\ref{Eq:UM3a}), together with the properties of $\sigma$, implies
that
\[
\big|X(t) - X(s)\big| \le \sigma\big(2^{-n}\big)\, n^{\frac{1 + \eps} \ga}
\le K_1\, \sigma\big(\rho(s, t)\big)\,\big(\log (1/\rho(s, t))\big)^{\frac{1 + \eps} \ga},
\quad \hbox{a.s.}
\]
for all $s, t \in [0, 1]^N$ such that $\rho(s, t)$ is small. Since $\eps > 0$ is arbitrary,
this proves (\ref{Eq:UM1d}).
\end{proof}

\section{Maximal moment index}
\label{Sec:MI}

It is clear that condition (\ref{Eq:MomCon1B}) [or (\ref{Eq:MomCon1})] is essential
for Theorem \ref{Thm:MC1}.
In many cases such as when $X =\{X(t), t \in T\}$ is a Gaussian or stable random field,
or has stationary increments, one can normalize the random variables
$X(\tau_n) - X(\tau_{n-1}')$ so that (\ref{Eq:MomCon1}) is reduced to conditions
on the maximal $\gamma$-moments of a sequence of ``homogeneous'' random variables.

Motivated by this, we introduce the notion of maximal $\gamma$-moment index for a
(not necessarily stationary) sequence of random variables which, in turn, provides
some sufficient conditions for (\ref{Eq:MomCon1}) to hold. 

\begin{definition}\label{Def:maxindex}
Let $\{\xi_k, k \ge 1\}$ be a sequence of random variables such that, for
some positive constants $\gamma $ and $K_\gamma$, $\E\big(|\xi_k|^\gamma\big) = K_\gamma$
for all $k \ge 1$. Let $M_n(\gamma) = \E\big(\max_{1 \le k \le n} |\xi_k|^\gamma \big)$.
Then the maximal $\gamma$-moment (upper) index of $\{\xi_k, k \ge 1\}$ is defined by
\begin{equation}\label{Def:MI}
\theta_\gamma = \limsup_{n \to \infty} \frac{\log M_n(\gamma)} {\log n}.
\end{equation}
If $\gamma = 1$, we simply call $\theta_1$ the maximal moment index of
$\{\xi_k, k \ge 1\}$.
\end{definition}

Some remarks are in order.
\begin{itemize}
\item If $0 < \gamma < \beta$, then Jensen's inequality implies that $\theta_\gamma
\le \frac{\ga} {\beta}\, \theta_\beta$.
\item Since $ \E(|\xi_1|^\gamma) \le M_n(\gamma) \le \E(\sum_{k=1}^n\, |\xi_k|^\gamma)$,
we clearly have  $\theta_\gamma \in [0, 1]$.

\item The maximal $\gamma$-moment
index carries some information about the dependence structure and distributional
properties of $\{\xi_k, k \ge 1\}$. Usually the choice of $\gamma$ is
determined by the heaviness of the tail probabilities of
$\{\xi_k, k \ge 1\}$. For a fixed $\gamma$, smaller
value of $\theta_\gamma$ indicates more dependence among $\{\xi_k, k \ge 1\}$.
This can be seen through the two extreme examples of stationary processes
$\{\xi_k, k \ge 1\}$: If $\xi_k \equiv \xi$, then $\theta_\gamma= 0$; while if
$\xi_k\ (k \ge 1)$ are i.i.d. and satisfy the conditions (\ref{Eq:TailP}) and
(\ref{Eq:TailP2}) below, then $\theta_\gamma= \gamma/\alpha$ for
$\gamma \in (0, \alpha)$.
Further evidence can be found in Samorodnitsky (2004) where it is shown that
the partial maxima of long-range dependent stable processes grow slower than
those of short-range dependent processes.
\end{itemize}

The following consequence of Theorem \ref{Thm:MC1} shows the usefulness of maximal
$\gamma$-moment index in determining the uniform modulus of continuity.

\begin{corollary}\label{Co:Ext-Kol}
Let $X = \{X(t), t \in \R^N\}$ be a continuous random field with values in $\R$ and
let $\{D_n, n \ge 1\} \subseteq [0, 1]^N$ be defined by
\begin{equation}\label{Def:Dn}
D_n = \left\{\Big(\frac{k_1} {2^n}, \cdots, \frac{k_N} {2^n}\Big):\
1 \le k_j \le 2^{n}, 1 \le j \le N\right\}.
\end{equation}
We assume there exist constants $\gamma\in (0, 1]$ and $H > 0$ such that for
all $s, t \in [0, 1]^N$
\begin{equation}\label{Eq:Ext-Kol1}
\E \big(|X(s) - X(t)|^\ga \big) \le  K\, |s-t|^{H\ga}.
\end{equation}
Consider the normalized random variables
\begin{equation}\label{Eq:Ext-Kol2}
\frac{X(\tau_p) - X(\tau_{p-1}') } {\big[\E (|X(\tau_p) - X(\tau'_{p-1})|^\ga)\big]^{1/\ga}},
\qquad \forall\ \tau_p \in D_p, \  \tau'_{p-1} \in O_{p-1}(\tau_{p})\
\hbox{ and }\ \forall\ p \ge 1,
\end{equation}
and number them according to the order $D_1$, $D_2\backslash D_1, \ldots$ and
denote the sequence by $\{\xi_k, k \ge 1\}$.  If $\{\xi_k, k \ge 1\}$  has a maximal
$\ga$-moment index $\theta:= \theta_\ga$ and $H \ga > N\theta$, then for every $\ep > 0$,
\begin{equation}\label{Coro:UM}
\limsup_{h\to 0+} \frac{\sup_{t \in [0,1]^N}\sup_{ |s-t|\le h}\big|X(t) - X(s)\big|}
{h^{H- \frac{N\theta} {\ga} - \ep }} = 0,
\qquad \hbox{ a.s.}
\end{equation}
Namely, $X(t)$ is uniformly H\"older continuous on $[0, 1]^N$ of all orders
$< H- \frac{N \theta} {\ga}$.
\end{corollary}

\begin{proof}\ For any $\ep > 0$, it follows from  (\ref{Def:MI}) that
$$M_{2^{Nn}}(\gamma) = \E\Big(\max_{1 \le k \le \#D_n} |\xi_k|^\gamma\Big)
 \le 2^{n (N \theta + \gamma \ep)}
 $$
for all $n$ large enough.  Combining this with (\ref{Eq:Ext-Kol1}) we derive
\[
\E\bigg(\max_{\tau_n \in D_n} \max_{\tau'_{n-1} \in O_{n-1}(\tau_{n}) }
\big|X(\tau_n) - X(\tau_{n-1}')\big|^\gamma\bigg)
\le K\, 2^{- (H - \frac{N\theta} \ga -\ep) \gamma n}
\]
for all integers $n$ large enough. Hence \eqref{Coro:UM} follows from Corollary
\ref{Thm:MC1b}.
\end{proof}

As an example, we show that the following multiparameter version of
Kolmogorov's continuity theorem follows from Corollary
\ref{Co:Ext-Kol}.

\begin{corollary}\label{Co:Kol}{\rm [{\bf Kolmogorov's continuity theorem}]}\
Let $X= \{X(t), t \in [0, 1]^N\}$ be a real-valued stochastic process.
Suppose there exist positive constants $\beta, K$ and $ \delta$ such that
\begin{equation}\label{Eq:Kol3}
\E \big(|X(s) - X(t)|^\beta \big) \le  K\, |s-t|^{N + \delta}, \qquad \forall \
s, t \in [0, 1]^N.
\end{equation}
Then $X$ has a version  which is uniformly H\"older
conditions on $[0, 1]^N$ of all orders $< \frac{\delta} {\beta}$.
\end{corollary}

\begin{proof}\ For any integer $n \ge 1$, let
$D_n$ be defined by (\ref{Def:Dn}) and let $\{\xi_k, k \ge 1\}$ be
the sequence of the normalized random variables in (\ref{Eq:Ext-Kol2}).
By (\ref{Eq:Kol3}), we see that (\ref{Eq:Ext-Kol1}) is satisfied with
$\gamma = \beta$ and $H = \frac{N+\delta}
{\beta}$.  In order to verify the second condition in Corollary \ref{Co:Ext-Kol},
we note that $\E\big(|\xi_k|^\ga\big) =1$ for every $k \ge 1$ and
use the trivial bound $\E\big(\max_{1 \le k \le m} |\xi_k|^\ga\big)\le m$
for all integers $m \ge 1$ to derive
$\theta_\ga\le 1$. Then it is clear that the conclusion of Corollary \ref{Co:Kol} follows
from Corollary \ref{Co:Ext-Kol}.
\end{proof}

In the rest of this section, we study the maximal moment indices for several
classes of random variables. The following lemma on Gaussian sequences will be
useful in the next section.

\begin{lemma}\label{Lem:GaussianMI}
Let $\{\xi_k, k \ge 1\}$ be a sequence of jointly Gaussian random variables
with mean 0 and variance 1. Then the following statements hold:
\begin{itemize}
\item[(i)]\ There is a universal constant $K_4>0$ such that
\begin{equation}\label{Eq:Gauss1}
\E\Big(\max_{1 \le k \le n} |\xi_k| \Big) \le K_4\, \sqrt{\log n}.
\end{equation}
In particular, for every $\gamma \in (0, 1]$, the maximal $\ga$-moment index of
$\{\xi_k, k \ge 1\}$ is 0.
\item[(ii)]\ If $|\E(\xi_j \xi_k)| \le \delta$ for a constant $\delta \in (0, 1)$
and for  all $1\le j < k\le n$, then there exists a constant $K_5$ such that
$\E\big(\max_{1 \le k \le n} |\xi_k| \big)
\ge K_5\, \sqrt{\log n}$.
\end{itemize}
\end{lemma}
\begin{proof}\ This lemma is well known, and we include a proof for completeness.
Part (i) can be proved by using the metric entropy method.
For any fixed integer $n \ge 1$, let $T = \{1, \ldots, n\}$ equipped with the
canonical metric $d (i, j) = \big[\E(\xi_i - \xi_j)^2\big]^{1/2}$. Then the
$d$-diameter of $T$ is at most 2. For any $\ep \in (0, 1)$, the $\ep$-covering
number $N_d(T, \ep) \le n$. Hence Dudley's entropy theorem [cf. e.g., Marcus and
Rosen (2006, Theorem 6.1.2) or Adler and Taylor (2007, Theorem 1.3.3)] gives
\begin{equation}\label{Eq:Gauss2}
\E\Big(\max_{1 \le k \le n} |\xi_k| \Big) \le K\, \int_0^1 \sqrt{\log n}\ d\ep,
\end{equation}
which yields \eqref{Eq:Gauss1}.

Under the condition of (ii), we have $d(i, j) = 2(1 - \E(\xi_i \xi_j))\ge 2 (1-\delta)$
for all $i \ne j$. Hence the conclusion of (ii) follows from the Sudakov minoration
[see Talagrand (2006, Lemma 2.1.2)].
\end{proof}

Dudley's entropy theorem has been extended to stochastic processes in
Orlicz spaces, see Talagrand (2006, p.30) for related references. Let
$\Psi: \R_+ \to \R_+$ be a convex function  such that $\Psi(0) = 0$
and $\Psi(r) > 0$ if $r \ne 0$. For a random variable $\xi$, the Orlicz
norm of $\xi$ is defined as
$$
\|\xi\|_\Psi = \inf\left\{c > 0: \E\Psi\big(\frac {\xi} c\big) \le 1\right\}.
$$
If $\{\xi_k, k \ge 1\}$ is a sequence of random variables such that
$\|\xi_k\|_\Psi \equiv K < \infty$, then (1.53) in Talagrand (2006) implies that
$\E\big(\max_{1 \le k \le n} |\xi_k| \big) \le K\,\Psi^{-1}(n) $.
This result extends Part (i) of Lemma \ref{Lem:GaussianMI}. By taking $\Psi(r)
= r^p$ ($p \ge 1$), one gets an upper bound for $\E\big(\max_{1 \le k \le n}
|\xi_k| \big)$ when $\{\xi_k, k \ge 1\}$ is a sequence of random variables
with the same finite $p$-th absolute moments.

The following lemma is applicable to $\a$-stable random variables.

\begin{lemma}\label{Lem:HeavyMI}
Let $\{\xi_k, k \ge 1\}$ be a sequence of random variables. 
The following statements hold:
\begin{itemize}
\item[(i)]\ If there exist positive constants $\alpha $ and $K_6$ such that
\begin{equation}\label{Eq:TailP}
\P\big( |\xi_k|\ge u\big) \le K_6\, u^{- \alpha}, \qquad \ \ \forall \ k \ge 1 \
\hbox{ and }\ u > 0,
\end{equation}
then for any $\gamma \in (0, \a)$ and $\ep > 0$ there is a finite constant $K_7$ such that
\begin{equation}\label{Eq:HeavyMI2}
\E\Big(\max_{1 \le k \le n} |\xi_k|^\ga \Big) \le K_7\, n^{\ga/\a} \,(\log n)^{(1+\ep)\ga/\a}
\end{equation}
for all integers $n \ge 2$. Consequently, for any $\gamma < \alpha$, we
have $\theta_\gamma \le \gamma/\alpha$.
\item[(ii)]\ If there exists a positive constant $K_8$
such that
\begin{equation}\label{Eq:TailP2}
\P\big( |\xi_k|\ge u\big) \ge K_8\, u^{- \alpha}, \qquad \ \ \forall \ k \ge 1\
\hbox{ and } \ u > 0
\end{equation}
and for all integers $n \ge 1$ and $u > 0$
\begin{equation}\label{Eq:FPOD}
\P\Big( \max_{1 \le k \le n}|\xi_k|\le u\Big) \le
\prod_{k=1}^n \P\big(|\xi_k|\le u\big),
\end{equation}
then, for any $\gamma \in (0, \a)$, there is a constant $K_9 > 0$ such that
\begin{equation}\label{Eq:HeavyMI3}
\E\Big(\max_{1 \le k \le n} |\xi_k|^\ga \Big) \ge K_9\, n^{\ga/\a}
\end{equation}
for all integers $n \ge 1$.
\end{itemize}
\end{lemma}

\begin{remark}
Clearly (\ref{Eq:FPOD}) is satisfied if the random variables $\xi_k\ (k \ge 1)$
are independent or, more generally, if the random variables $|\xi_k|\ (k \ge 1)$
are negatively orthand dependent, that is, for all $n \ge 1$ and $x_1, \ldots,
x_n >0$,
\[
\P\Big(|\xi_1| \le x_1, \ldots, |\xi_n| \le x_n\Big)
\le \prod_{k=1}^n \P\big(|\xi_k| \le x_k
\big).
\]
Moreover, by modifying the proof of Part (ii) of Lemma \ref{Lem:HeavyMI},
one can show that, if we assume in (i)
that the random variables $|\xi_k|\ (k \ge 1)$ satisfies
$$
\P\Big( \max_{1 \le k \le n}|\xi_k|\le u\Big) \ge
\prod_{k=1}^n \P\big(|\xi_k|\le u\big),
$$
then (\ref{Eq:HeavyMI2}) can be improved to $\E\big(\max_{1 \le k \le n}
 |\xi_k|^\ga \big) \le K\, n^{\ga/\a}.$
\end{remark}

\begin{proof}{\bf of Lemma \ref{Lem:HeavyMI}}\ Let $\gamma \in (0, \alpha)$ and $\ep > 0$
be given constants. To prove Part (i),
let $\Phi: \R_+ \to \R_+$ be a nondecreasing and convex function
satisfying $\Phi(0) = 0$ and
\begin{equation}\label{Eq:phi1}
\Phi(x) \sim \frac{x^{\a/\ga}} {(\log x)^{1+\ep}} \qquad \hbox{ as }\ x \to \infty.
\end{equation}
Here and in the sequel, $f(x) \sim g(x)$ means $f(x)/g(x) \to 1$
as $x \to \infty$ or $x \to 0$. Then the inverse function of $\Phi(x)$, denoted by
$\Phi^{-1}(x)$, is nonnegative, nondecreasing and concave on $[0, \infty)$. Moreover,
\begin{equation}\label{Eq:phi2}
\Phi^{-1}(x) \sim \frac \ga \a\, x^{\ga/\a} \,(\log x)^{(1+\ep) \ga/\a}
\qquad \hbox{ as }\ x \to \infty.
\end{equation}

By (\ref{Eq:TailP}) and (\ref{Eq:phi1}) we derive $\E\big(\Phi(|\xi_k|^\ga)\big)
\le K_{10}$ for all $k \ge 1$, where $K_{10}$ is a positive constant depending
on $K_6$, $\a$, $\ep$ and $\ga$ only. This and Jensen's inequality together imply
\[
\begin{split}
\E\Big(\max_{1 \le k \le n} |\xi_k|^\ga \Big) &\le
\Phi^{-1}\bigg[\E \Phi\Big(\max_{1 \le k \le n} |\xi_k|^\ga \Big)\bigg]\\
&\le \Phi^{-1}\bigg[ \sum_{1 \le k \le n}\E\big(\Phi(|\xi_k|^\ga)\big)\bigg]
= \Phi^{-1}(K_{10}\,n).
\end{split}
\]
Combining this and (\ref{Eq:phi2}) yields (\ref{Eq:HeavyMI2}).

To prove Part (ii), we write
\begin{equation}\label{Eq:int-by-parts}
\E\Big(\max_{1 \le k \le n} |\xi_k|^\ga \Big)= \ga \int_0^\infty u^{\ga - 1}
\P\Big(\max_{1 \le k \le n} |\xi_k| > u\Big)\, du.
\end{equation}
It follows from (\ref{Eq:FPOD}) and (\ref{Eq:TailP2}) that
\begin{equation}\label{Eq:int-by-parts-2}
\begin{split}
\P\Big(\max_{1 \le k \le n} |\xi_k| > u\Big) &\ge 1 -  \prod_{k=1}^n
\Big( 1 - \P\big(|\xi_k| > u\big)\Big) \\
&\ge 1 - \Big( 1- K_8\,u^{-\a}\Big)^n
\end{split}
\end{equation}
for all $u> K_8^{1/\alpha}$. By using the elementary inequality
$$1 - (1-x)^n \ge \frac 1 2 nx \qquad  \forall\ 0 \le x \le
1 - (\frac 1 2 )^{1/(n-1)},
$$
we derive from (\ref{Eq:int-by-parts-2}) that
\begin{equation}\label{Eq:int-by-parts-3}
\P\Big(\max_{1 \le k \le n} |\xi_k| > u\Big) \ge \frac 1 2\, n u^{-\a}
\end{equation}
for all $u \ge K\,\big(1 - (\frac 1 2 )^{1/(n-1)}\big)^{-1/\a} \asymp
n^{1/\a}$. Combining (\ref{Eq:int-by-parts}),
(\ref{Eq:int-by-parts-2})
and (\ref{Eq:int-by-parts-3}) we obtain
\begin{equation}\label{Eq:int-by-parts-4}
\E\bigg(\max_{1 \le k \le n} |\xi_k|^\ga \bigg)\ge K\,n\,
\int_{n^{1/\a}}^\infty u^{\ga - \a- 1}
\, du = K\, n^{\ga/\a}.
\end{equation}
This finishes the proof of Lemma \ref{Lem:HeavyMI}.
\end{proof}

Part (ii) of Lemma \ref{Lem:HeavyMI} indicates that, in general, it is
difficult to determine an optimal lower bound for $\E\big(\max_{1 \le
k \le n} |\xi_k|^\ga \big)$ without additional information
about the dependence structure of $\{\xi_k, k \ge 1\}$. We point out that,
even for stationary stable sequences, it is an open problem in general
to determine sharp (i.e., up to constant factors) upper and lower bounds for
$\E\big(\max_{1 \le k \le n} |\xi_k|^\ga \big)$.

Here are some partial results. For symmetric stable
or infinitely divisible sequences, some lower bounds can
be derived from Theorem 2.3.1 and Theorem 5.3.2 in Talagrand (2006). For example,
if $\{\xi_k\}$ is a stationary sequence of symmetric $\a$-stable (S$\alpha$S)
random variables such that $\a \in (1, 2]$ and there exists a constant $\delta > 0$ such that
$\|\xi_k - \xi_j\|_\a \ge \delta$ for all $k \ne j$ (where
$\|\xi\|_\a$ denotes the scale parameter of $\xi$), then
$\E\big(\max_{1 \le k \le n} |\xi_k|\big) \ge K (\log n)^{(\a -1)/\a}$. On the
other hand, Samorodnitsky (2004) studied the rate of growth of the partial maxima
$M_n = \max\{|\xi_k|: 1 \le k \le n\}$ of a stationary $\alpha$-stable
sequence $\{\xi_k, k \ge 1\}$ based on the ergodic theoretical properties of the
underlying flow. He discovered that (i) if the stationary S$\alpha$S process $\{\xi_k, k \ge 1\}$
is generated by a dissipative flow then $M_n$ grows always at the rate of $n^{1/\alpha}$ and (ii)
if the stationary S$\alpha$S process $\{\xi_k, k \ge 1\}$
is generated by a conservative flow then $M_n$ grows at the rate slower than $n^{1/\alpha}$.
Samorodnitsky (2004, p.1440) conjectured that many other important
properties of $\{\xi_k, k \ge 1\}$ will also change as the underlying flow changes from
being dissipative to being conservative. In particular, we believe that, if  $X= \{X(t), t \in \R\}$
is a (self-similar) S$\alpha$S process with stationary increments, then the uniform
modulus of continuity of $X$ depends on the nature of the flow generating the stationary
sequence $\xi_k = X(k+1) - X(k)$ $(k \ge 0)$. This idea will be pursued further elsewhere.

Combining Lemma \ref{Lem:HeavyMI} with the proof of Theorem \ref{Thm:MC1},
we have the following result on modulus of continuity for general self-similar
processes with stationary increments, which is an improvement of
Corollary \ref{Co:Ext-Kol}. In the special case of $N=1$ and
$X = \{X(t), t \in \R\}$ being an $\alpha$-stable process, it recovers
Theorem 2 in K\^ono and Maejima (1991b).
\begin{corollary}\label{Co:Kol2}
Let $X= \{X(t), t \in \R^N\}$ be a real-valued $H$-self-similar random field
with stationary increments. Let $V= \{v_\ell, 1 \le \ell \le 2^N-1\}$ be the
set of vertices of $[0, 1]^N$, excluding $0$. If there exists a constant
$\alpha > \frac N H$ such that
\begin{equation}\label{Eq:TailP3}
\P\big( |X(v_\ell)|\ge u\big) \le K\, u^{- \alpha}, \qquad \forall \, v_\ell \in V \,
\hbox{ and } \ u > 0,
\end{equation}
then for any $\ep > 0$,
\begin{equation}\label{Coro:UM2}
\limsup_{h\to 0+} \frac{\sup_{t \in [0,1]^N}\sup_{ |s-t|\le h}\big|X(t) - X(s)\big|}
{h^{H - \frac N \a} \, \big(\log 1/h\big)^{\frac{2+\ep} \a}} = 0,
\qquad \hbox{ a.s.}
\end{equation}
\end{corollary}

\begin{proof}\ Again we use the $\ell^\infty$ metric in $\R^N$. For every $n \ge 1$,
let $D_n$ be defined as in (\ref{Def:Dn}). Then the sequence $\{D_n, n \ge 1\}$
satisfies the conditions in Section 2. The
self-similarity and stationarity of the increments of $X$ imply that, for every
$\tau_n \in D_n$
and $\tau_{n-1}' \in O_{n-1}(\tau_n)$, there is a $v_\ell \in V$ such that
\[
X(\tau_n ) - X(\tau_{n-1}' )
\stackrel{d}{=} |\tau_n - \tau_{n-1}'|^H X(v_\ell).
\]
Hence we  apply Lemma \ref{Lem:HeavyMI} to the normalized sequence
$\{\frac{X(\tau_n ) - X(\tau_{n-1}' )}
{|\tau_n - \tau_{n-1}'|^H}\}$ to derive that for any $0 < \ga < \min\{1, \a\}$
\begin{equation}\label{Eq:SSSI-1}
\E\bigg(\max_{\tau_n \in D_n} \max_{\tau'_{n-1} \in O_{n-1}(\tau_{n})}
\big|X(\tau_n) - X(\tau_{n-1}')\big|^\gamma\bigg)
\le K\, 2^{-(H - \frac N \alpha) \gamma n} \, \big(\log 2^{Nn}\big)^{(1+\ep)\ga/\a}.
\end{equation}
Hence (\ref{Coro:UM2}) follows from (\ref{Eq:SSSI-1}) and Theorem \ref{Thm:MC1}.
\end{proof}

\section{Applications to stable random fields}
\label{Sec:Appl}

In this section  we apply the results in Sections \ref{sec:MC} and \ref{Sec:MI}
to harmonizable-type $\a$-stable random fields with $0 < \a < 2$. Similar methods
can be applied to other types of stable random fields, or more generally,
infinitely divisible processes. 

\subsection{Harmonizable-type stable random fields with stationary increments}

Let $X=\{X(t),\, t \in \R^N\}$ be a real-valued stable random
field defined by
\begin{equation}\label{StableF2}
X(t) = \Re \int_{\R^N} \big(e^{i\l t, x\r} - 1\big)\, \widetilde{M}_\alpha(dx),
\end{equation}
where $\widetilde{M}_\alpha$ is a rotationally invariant $\a$-stable random
measure on $\R^N$ with control measure $\Delta$, which satisfies
\begin{equation}\label{Eq:LM}
\int_{\R^N} (1 \wedge |x|^\a)\, \De(d x) < \infty.
\end{equation}
This condition assures that stochastic integral in (\ref{StableF2})
is well-defined, see Samorodnitsky and Taqqu (1994, Chapter 6) for
further information. The measure $\De$ is called the spectral measure
of $X$ and its density function, when it exists, is called the spectral
density of $X$.

It can be verified that the stable random field $X$ defined by
(\ref{StableF2}) has stationary increments and $X(0) = 0$.
Denote the scale parameter of $X(t)$ by $\|X(t)\|_\a$. Then for all
$t \in \R^N$,
\begin{equation}\label{Eq:ScaleP}
\|X(t)\|_\a^\a = 2^{\a/2}
\int_{\R^N} \bigl(1 - \cos \l t, x \r \bigr)^{\a/2} \, \De(d x).
\end{equation}
Similar to Gaussian processes, this function plays an
important role in studying sample path properties of stable random field
$X$ defined by (\ref{StableF2}).

For simplicity, we will assume that the spectral measure
$\De$ is absolutely continuous and its density function $f(x)$
satisfies the following condition
\begin{equation}\label{Eq:spden1}
f(x) \le K_{11}\, |x|^{-(\a H + N)}, \qquad \quad \forall \, x \in \R^N\
\hbox{ with }\ |x| \ge K_{12},
\end{equation}
where $K_{11}, K_{12} > 0$ and $H \in (0, 1)$ are constants. As shown by Theorem
\ref{Thm:HarmF} below, the parameter $H$ determines the smoothness
of the sample function $X(t)$.

Now we provide some examples of stable random fields
satisfying the condition (\ref{Eq:spden1}).

\begin{example}\ [Harmonizable fractional stable motion]\
Let $H\in (0, 1)$ and $\a\in (0, 2]$ be given constants.
The harmonizable fractional stable field $\widetilde Z^{H} = \{
\widetilde Z^{H}(t), t \in \R^N\}$ with values in $\R$ is
defined by (\ref{StableF2}) with spectral density
\begin{equation}\label{Eq:fsm1}
f_{H, \alpha}(x) = c(\a,H, N) \, \frac 1 {|x|^{\a H + {N} }},
\end{equation}
where $c(\a,H, N)>0$ is a normalizing constant such that the
scale parameter of $\widetilde Z^{H}(e_1)$ equals 1, where
$e_1=(1, 0, \ldots, 0)\in \R^N$. Hence we have $\|\widetilde Z^{H}(t)\|_\a
= |t|^{H}$ for all $t \in \R^N$.

It is easy to verify that the $\a$-stable random field $\widetilde
Z^{H}$ is $H$-self-similar with stationary and isotropic increments
[or \emph{stationary increments in the strong sense} in terms of
Samorodnitsky and Taqqu (1994, p.392)]. The random field $\widetilde
Z^{H}$ is a stable analogue of fractional Brownian motion $B^H$ of
index $H$ and serves as an important representative for understanding
harmonizable-type stable random fields. Even though the results
on local times of $\widetilde Z^{H}$ [Nolan (1989) and Xiao (2008)]
and uniform modulus of continuity [see Theorem \ref{Thm:HarmF} below]
show that $\widetilde Z^{H}$ shares many sample path properties with
$B^H$, little has been known about the sharpness of these results and
the existing tools do not seem to be capable for attacking these
problems. It is an interesting (and challenging) task to develop
new methods for studying the fine structures of $\widetilde Z^{H}$
and other stable random fields.
\end{example}

\begin{example} \ [Fractional Riesz-Bessel $\alpha$-stable motion]\
Consider the stable random field $X= \{X(t), t \in
\R^N\}$ in $\R$ defined by (\ref{StableF2}) with spectral density
\begin{equation}\label{Eq:Anhdensity}
f_{\gamma, \eta}(x) = \frac {c(\a, \ga, \eta, N)} {|x|^{2 \ga} (1
+ |x|^2)^\eta},
\end{equation}
where  $\eta$ and $\gamma$ are positive constants satisfying
\[
\eta + \ga > \frac{N} 2, \quad 0 < 2 \ga < \a + N
\]
and $c(\a, \ga, \eta, N)>0$ is a normalizing constant. When $\a = 2$
[i.e., the Gaussian case]
such density functions were consider by Anh et al. (1999). Since
$f_{\gamma, \eta}$ involves both the Fourier transforms of the Riesz
kernel and the Bessel kernel, Anh et al. (1999) called the corresponding
Gaussian random field $X$ the fractional Riesz-Bessel
motion with indices $\eta $ and $\gamma$. They showed that
these Gaussian random fields can be used for modeling simultaneously
long range dependence and intermittency.

In analogy to the terminology in  Anh et al. (1999), we call
$X$ the fractional Riesz-Bessel $\alpha$-stable motion with
indices $\eta $ and $\gamma$. Clearly, the spectral density
$f_{\gamma, \eta}(x)$ in (\ref{Eq:Anhdensity}) satisfies (\ref{Eq:spden1})
with $H = (2(\eta + \gamma) - N)/\a$. Moreover, since the spectral density
$f_{\gamma, \eta}(x)$ is regularly varying at infinity of order
$2 (\eta + \gamma) > N$, by modifying the proof of Theorem 1 in
Pitman \cite{Pitman68} we can show that, if $2(\gamma +
\eta) - N < \a$,  then $\|X(t)\|_\a$ is regularly varying at
0 of order $(2(\eta + \gamma) - N)/\a$ and
\begin{equation}\label{Eq:Pitman}
\| X(t)\|_\a \sim  \, |t|^{(2(\eta + \gamma) - N)/\a}, \qquad
\hbox{ as }\ |t| \to 0.
\end{equation}
We will see that the modulus of continuity of $X$ is determined by
(\ref{Eq:Pitman}).
\end{example}

The following result provides information on the maximal moment
of harmonizable-type stable random fields. It is more precise than
Part (i) of Lemma \ref{Lem:HeavyMI}.

\begin{proposition}\label{Prop:HarmF1}
Let $X = \{X(t), t \in [0, 1]^N\}$ be an $\alpha$-stable random field defined
by (\ref{StableF2}) with spectral density satisfying (\ref{Eq:spden1}).
Then for every $0 < \ga < \min\{1, \a\}$, $\eta > 0$ and $n \ge 1$,
\begin{equation}\label{Eq:HSF-1}
\E\bigg(\max_{\tau_n \in D_n} \max_{\tau'_{n-1} \in O_{n-1}(\tau_{n})}
\big|X(\tau_n) - X(\tau_{n-1}')\big|^\gamma\bigg)\le K\, 2^{-H\gamma n}\,
n^{\frac {(1+\eta)\ga} \a},
\end{equation}
where $\{D_n, n \ge 1\}$ is the sequence defined in (\ref{Def:Dn}).
\end{proposition}

In order to prove Proposition \ref{Prop:HarmF1}, we will make use of a
LePage-type representation for $X$, which allows us to view $X$ as a mixture of
Gaussian process. This powerful idea was due to Marcus and Pisier (1984) and has
become a standard tool for studying sample path regularity of stable processes.
For more information on series representations of infinitely divisible processes
and their applications, see Rosi\'nski (1989, 1990), K\^ono and Maejima (1991a),
Samorodnitsky and Taqqu (1994), Marcus and Rosinski (2005), Bierm\'e and Lacaux (2009),
just to mention a few.

We need some notation. Let $\mu$ be an arbitrary probability on $\R^N$ which
is equivalent to the Lebesgue measure $\lambda_N$. Denote by $\varphi=
d\mu/d\lambda_N$ the Radon-Nikodym derivative of $\mu$, so $\mu(dx) =
\varphi(x) \, dx$.
Assume that
\begin{itemize}
\item $\{\Gamma_j, j \ge 1\}$ is a sequence of Poisson
arrival times with intensity 1;
\item $\{g_j, j \ge 1\}$ is a sequence of i.i.d. complex
valued Gaussian random variables such that $g_j \stackrel{d} {=} e^{i \theta} g_j$
for all $\theta \in \R$ and $\E(|\Re g_1|^\a) =1$;
\item $\{\xi_j, j \ge 1\}$ is a sequence of  i.i.d. random variables with values in $\R^N$
and density function $\varphi$;
\item the sequences $\{\Gamma_j, j \ge 1\}$, $\{g_j, j \ge 1\}$ and $\{\xi_j, j \ge 1\}$
are independent. We denote the expectations with respect to $\{\Gamma_j, j \ge 1\}$,
$\{g_j, j \ge 1\}$ and $\{\xi_j, j \ge 1\}$ by $\E_\Gamma$, $\E_g$ and $\E_\xi$,
respectively.
\end{itemize}

The following lemma is from Bierm\'e and Lacaux (2009); see also K\^ono and
Maejima (1991a) and Marcus and Pisier (1984).

\begin{lemma}\label{Lem:Lepage}
For any family of complex-valued functions $h(t, \cdot) \in L^\a (\R^N, dx) \,
(0 < \alpha < 2)$, let $Z = \{Z(t), t \in \R^N\}$ be the $\a$-stable random
field defined by
\[
Z(t) = \Re \int_{\R^N} h(t, x)\, \widetilde{Z}_\a(dx), \qquad \forall t \in \R^N,
\]
where $\widetilde{Z}_\a$ is a complex-valued, rotationally invariant
$\alpha$-stable random measure on $(\R^N, {\cal B}(\R^N))$
with Lebesgue control measure. Then
\[
\bigg\{Z(t), t \in \R^N\bigg\}\,
\stackrel{d} {=}\, \bigg\{Y(t), t \in \R^N\bigg\},
\]
where $\stackrel{d} {=}$ means equality in finite dimensional distributions and
\begin{equation}\label{Eq:LePage2}
Y(t) = C_\a\, \Re\bigg(\sum_{j=1}^\infty \Gamma_j^{-1/\alpha}
\varphi(\xi_j)^{-1/\alpha}\, h(t, \xi_j)\, g_j\bigg).
\end{equation}
In the above, for every $t \in \R^N$, the random series (\ref{Eq:LePage2})
converges almost surely and $C_\a$ is the constant given by
\[
C_\a = \bigg(\frac 1 {2\pi} \int_0^{\pi} \big|\cos \theta \big|^\a\, d\theta\bigg)^{1/\a}
\bigg(\int_0^\infty \frac{\sin \theta} {\theta^\a}\, d\theta\bigg)^{-1/\a}.
\]
\end{lemma}

\begin{proof}{\bf of Proposition \ref{Prop:HarmF1}}\  Choose two positive constants
$\eta $ and $\beta$ such that $N-\frac{2 \a (1-H)} {2 -\a} < \beta < N$. Let $\varphi:
\R^N \backslash\{0\}
\to [0, \infty)$ be the function defined by
\begin{equation}\label{Eq:phi-1}
\varphi(x) = \left\{\begin{array}{ll}
K_{13}\, |x|^{-\beta}\qquad \qquad \qquad \qquad &\hbox{ if }\ |x| \le 3,\\
K_{14}\, \big(|x|^N(\log |x|)^{1 + \eta}\big)^{-1}
&\hbox{ if }\ |x| > 3,\\
\end{array}
\right.
\end{equation}
where the constants $K_{13}$ and $K_{14}$ are chosen such that
$\int_{\R^N}\varphi(x) \, dx = 1$.

By Lemma \ref{Lem:Lepage}, the stable random field $X$ defined by (\ref{StableF2})
with spectral density $f(x)$ has the same finite dimensional distributions as
\begin{equation}\label{Eq:LePage3}
Y(t) = C_\a\, \Re\bigg(\sum_{j=1}^\infty \Gamma_j^{-1/\alpha}
\varphi(\xi_j)^{-1/\alpha}\, h(t, \xi_j)\, g_j\bigg),
\end{equation}
where the function $h(t, x)$ is defined by $h(t, x)
=\big(e^{i\l t, x\r} - 1\big) f(x)^{1/\a}.$ Hence it is sufficient to prove
(\ref{Eq:HSF-1}) for $Y$.

Conditional on $\{(\xi_j, \Gamma_j),\, j \ge 1\}$, $Y$ is a Gaussian random field with
incremental variance given by
\begin{equation}\label{Eq:IncVar1}
\begin{split}
\E_g \Big[\big(Y(t) - Y(s)\big)^2 \Big] &= 2\, C_\a^2\,
\sum_{j=1}^\infty \Gamma_j^{-2/\alpha}
\varphi(\xi_j)^{-2/\alpha}\, \big(1 - \cos \langle t-s, \xi_j\r\big)\,f(\xi_j)^{2/\a}\\
&\le K \,\sum_{j=1}^\infty \Gamma_j^{-2/\alpha}
\varphi(\xi_j)^{-2/\alpha}\, \big\{1 \wedge  |t-s|^2 |\xi_j|^2\big\}\,
|\xi_j|^{-2(H + \frac N \a)},
\end{split}
\end{equation}
where the inequality follows from (\ref{Eq:spden1}) and $K>0$ is a constant.

Hence by using (\ref{Eq:IncVar1}) and Lemma \ref{Lem:GaussianMI}, we have
\begin{equation}\label{Eq:IncVar2}
\begin{split}
&\E_g\bigg(\max_{\tau_n \in D_n} \max_{\tau'_{n-1} \in O_{n-1}(\tau_{n})}
\big|Y(\tau_n) - Y(\tau_{n-1}')\big|^\gamma\bigg)\\
&\le K\, \max_{\tau_n \in D_n} \max_{\tau'_{n-1} \in O_{n-1}(\tau_{n})}
 \Big[ \E_g\big(Y(\tau_n) - Y(\tau_{n-1}')\big)^2\Big]^{\gamma/2}\, n^{\ga/2}\\
&\le K\, \bigg[\sum_{j=1}^\infty \Gamma_j^{-2/\alpha}
\varphi(\xi_j)^{-2/\alpha}\, \big\{1 \wedge  2^{-2n} |\xi_j|^2\big\}\,
|\xi_j|^{-2(H + \frac N \a)} \bigg]^{\gamma/2}\, n^{\ga/2}.
\end{split}
\end{equation}
It remains to show that (\ref{Eq:HSF-1}) follows from taking expectations
on both sides of (\ref{Eq:IncVar2}) with respect to $\{\xi_j, j \ge 1\}$
and $\{\Gamma_j, j \ge 1\}$.

Note that, for every $j \ge 1$, $\Ga_j$ is a Gamma random variable
with density function
$$
p(x) = \frac{x^{j-1}\, e^{-x}} {(j-1)!}, \qquad \forall \ x \ge 0.
$$
It is elementary to verify that there is a constant $K>0$ such that
\begin{equation}\label{Eq:IncVar3a}
\E_\Ga\big(\Ga_j^{-2/\a}\big) \le
K\,j^{-2/\a},\quad \ \ \ \forall \ j > 2/\a
\end{equation}
and
\begin{equation}\label{Eq:IncVar3b}
\E_\Ga\big(\Ga_j^{-\ga/\a}\big) \le K \qquad \ \forall\ 1 \le j \le 2/\a.
\end{equation}

Since $\gamma \in (0, 1)$, we use the elementary inequality $(x +y)^{\ga/2}
\le x^{\ga/2} + y^{\ga/2}$ and
Jensen's inequality to derive
\begin{equation}\label{Eq:IncVar3}
\begin{split}
&\E_{\Gamma, \xi} \Bigg\{\bigg(\sum_{j=1}^\infty \Gamma_j^{-2/\alpha}
\varphi(\xi_j)^{-2/\alpha}\, \big\{1 \wedge  2^{-2n} |\xi_j|^2\big\}\,
|\xi_j|^{-2(H + \frac N \a)} \bigg)^{\gamma/2}\Bigg\} \\
&\le \E_{\Gamma, \xi} \bigg(\sum_{j=1}^{\lfloor 2/\a \rfloor} \Gamma_j^{-\ga/\alpha}
\varphi(\xi_j)^{-\ga/\alpha}\, \big\{1 \wedge  2^{-\ga n} |\xi_j|^\ga\big\}\,
|\xi_j|^{-\ga(H + \frac N \a)}\bigg) \\
&\qquad \quad +\bigg[ \E_{\Gamma, \xi} \bigg(\sum_{j=\lfloor 2/\a \rfloor +1}^\infty
\Gamma_j^{-2/\alpha}
\varphi(\xi_j)^{-2/\alpha}\, \big\{1 \wedge  2^{-2n} |\xi_j|^2\big\}\,
|\xi_j|^{-2(H + \frac N \a)}\bigg) \bigg]^{\gamma/2}\\
&:= I_1 + I_2.
\end{split}
\end{equation}

Since  $I_1$ and $I_2$ can be estimated by using the same method, we only consider
$I_2$ below. Taking the expectation $\E_\xi$ first, we have
\begin{equation}\label{Eq:IncVar4}
\begin{split}
&\E_\xi\bigg(\varphi(\xi_j)^{-2/\alpha}\, \big\{1 \wedge  2^{-2n} |\xi_j|^2\big\}\,
|\xi_j|^{-2(H + \frac N \a)}\bigg)\\
&= \int_{\R^N}\varphi(x)^{1 - \frac 2 \a} \, \big(1 \wedge  2^{-2n} |x|^2\big)\, \frac{dx}
{|x|^{2(H + \frac N \a)}}\\
&\ \ = 2^{-2n}\, \int_{|x| \le 2^n}  \varphi(x)^{1 - \frac 2 \a} \,\frac { dx}
{|x|^{2(H + \frac N \a -1)}} + \int_{|x| > 2^n}  \varphi(x)^{1 - \frac 2 \a}
\,\frac {dx} {|x|^{2(H + \frac N \a)}} \\
&\ \ := J_1 + J_2.
\end{split}
\end{equation}
By (\ref{Eq:phi-1}) and a change of variables, we derive
\begin{equation}\label{Eq:IncVar5}
\begin{split}
J_1 &= K \, 2^{-2n}\,\bigg(\int_0^3 \frac{dr} {r^{2H + (N-\beta) (\frac 2 \a-1)- 1}} \\
&\qquad \qquad \qquad \qquad +  \int_3^{2^n} \frac{dr} {r^{2H  -1}
|\log r|^{(1 + \eta)(1 - \frac 2 \a)}}\bigg)\\
&\le K\, 2^{-2Hn}\, n^{- (1 + \eta)(1 - \frac 2 \a)}.
\end{split}
\end{equation}
Note that, because of the choice of $\beta$, the first integral is constant.
Similarly, we also have
\begin{equation}\label{Eq:IncVar6}
\begin{split}
J_2 \le K \, 2^{-2Hn}\, n^{-(1 + \eta)(1 - \frac 2 \a)}.
\end{split}
\end{equation}
By (\ref{Eq:IncVar3a}), (\ref{Eq:IncVar4}), (\ref{Eq:IncVar5}) and (\ref{Eq:IncVar6}),
we obtain
\begin{equation}\label{Eq:IncVar7}
\begin{split}
I_2 &\le \bigg[\E_{\Gamma} \bigg(K\, \sum_{j= \lfloor 2/\a \rfloor + 1}^\infty
\Gamma_j^{-2/\alpha}  2^{-2Hn}\, n^{-(1 + \eta)(1 -\frac 2 \a)} \bigg) \bigg]^{\ga/2}\\
&\le K \, 2^{-\ga Hn}\, n^{-\ga (1 + \eta)(\frac 1 2 - \frac 1 \a)}\,
\bigg[\sum_{j = \lfloor 2/\a \rfloor + 1 }^\infty
\E_\Ga\Big(\Ga_j^{-2/\a}\Big) \bigg]^{\ga/2}\\
&\le K\, 2^{-\ga Hn}\, n^{-\ga (1 + \eta)(\frac 1 2 - \frac 1 \a)}.
\end{split}
\end{equation}
Similarly, we can derive
\begin{equation}\label{Eq:IncVar8}
I_1 \le K\,  2^{-\ga Hn}\, n^{-\ga (1 + \eta)(\frac 1 2 - \frac \ga \a)}.
\end{equation}
Combining (\ref{Eq:IncVar3}), (\ref{Eq:IncVar7}) and (\ref{Eq:IncVar8}) we obtain
\begin{equation}\label{Eq:IncVar9}
\begin{split}
&\E_{\Gamma, \xi} \Bigg\{\bigg(\sum_{j=1}^\infty \Gamma_j^{-2/\alpha}
\varphi(\xi_j)^{-2/\alpha}\, \big\{1 \wedge  2^{-2n} |\xi_j|^2\big\}\,
|\xi_j|^{-2(H + \frac N \a)} \bigg)^{\gamma/2}\Bigg\} \\
& \le K\,  2^{-\ga Hn}\, n^{-\ga (1 + \eta)(\frac 1 2 - \frac \ga \a)}.
\end{split}
\end{equation}
Finally  (\ref{Eq:HSF-1}) follows from (\ref{Eq:IncVar2}) and (\ref{Eq:IncVar9}).
This proves Proposition \ref{Prop:HarmF1}.
\end{proof}

The following is a consequence of Proposition \ref{Prop:HarmF1} and Theorem
\ref{Thm:MC1}.

\begin{theorem}\label{Thm:HarmF}
Let $X = \{X(t), t \in \R^N\}$ be an $\alpha$-stable random field defined by
(\ref{StableF2}) with spectral density satisfying (\ref{Eq:spden1}).
Then for any $\ep > 0$,
\begin{equation}\label{Eq:MC-Harm}
\limsup_{h\to 0+} \frac{\sup_{t \in [0,1]^N}\sup_{ |s-t|\le h}\big|X(t) - X(s)\big|}
{h^{H} \, \big(\log 1/h\big)^{(2 + \varepsilon)/\a} } = 0,
\qquad \hbox{ a.s.}
\end{equation}
\end{theorem}

\begin{proof}\ It follows from Proposition \ref{Prop:HarmF1} that for all $\eta>0$ and
integers $n \ge 1$,
\[
\sum_{p=n}^\infty \E\bigg(\max_{\tau_p \in D_p} \max_{\tau'_{p-1} \in O_{p-1}(\tau_{p})}
\big|X(\tau_p) - X(\tau_{p-1}')\big|^\gamma\bigg)
\le K\, 2^{- H\ga n}\, n^{\frac {\ga(1 +\eta)} \a}.
\]
Thus, $X$ satisfies (\ref{Eq:MomCon1B}) with  $\sigma(h) =
h^{H} (\log 1/h)^{\frac {1+\eta} \a}$. Since $\eta>0$ and $0 < \ga < \a$ are arbitrary,
(\ref{Eq:MC-Harm}) follows from Theorem \ref{Thm:MC1}.
\end{proof}

\subsection{Harmonizable fractional stable sheets}
\label{sec:anisotropic}

For any given $0 < \alpha < 2$ and  $\vec H = (H_1,$ $\ldots, H_N)\in (0,
1)^N$, we define the harmonizable fractional stable sheet
$\widetilde Z^{\vec H} = \{ \widetilde Z^{\vec H}(t), t \in
\R_+^N\}$ with values in $\R$ by
\begin{equation}\label{eq:stableharm2}
 \widetilde Z^{\vec H} (t) =  \Re \int_{\R^N} \prod_{j=1}^N\,
\frac{e^{it_j x_j}-1}{|x_j|^{H_j+\frac{1}{\a}}} \,
\widetilde{Z}_\a(d\lambda),
\end{equation}
where $\widetilde{Z}_\alpha$ is a complex-valued random
measure as in Lemma \ref{Lem:Lepage}.

From (\ref{eq:stableharm2}) it follows that $\widetilde Z^{\vec H}$
has the following operator-scaling property:
For any $N\times N$ diagonal matrix $E = (b_{ij})$ with
$b_{ii} = b_i > 0$ for all $1 \le i \le N$ and $b_{ij}=0$ if $i\ne
j$, we have
\begin{equation}\label{Eq:OSS}
\big\{ \widetilde Z^{ \vec H}(E t),\, t \in \R^N \big\}
\stackrel{d}{=} \bigg\{ \bigg(\prod_{j=1}^N b_j^{H_j}\bigg)
\widetilde Z^{\vec H}(t),\ t \in \R^N \bigg\}.
\end{equation}
Along each direction of $\R^N_+$, $\widetilde{Z}^{\vec H}$ becomes a
real-valued harmonizable fractional stable motion [cf. Samorodnitsky
and Taqqu (1994, Chapter 7)]. When the indices $H_1, \ldots, H_N$ are
not the same, $\widetilde Z^{\vec H}$ has different scaling behavior
along different directions and this anisotropic nature induces some
interesting geometric and analytic properties for $\widetilde Z^{\vec H}$.
Note that $\widetilde Z^{\vec H}$ does not have stationary increments in the
ordinary sense and, thus, is different from the operator-scaling
stable fields considered in Bierm\'e and Lacaux (2009).
See Xiao (2006, 2008) for further information on sample path properties of
stable random fields.

The following result gives the uniform modulus of continuity for
$\widetilde Z^{\vec H}$, which is significantly different from the
result for linear fractional stable sheets obtained in Ayache, Roueff and Xiao
(2007, 2009).

\begin{theorem} \label{thm:HFSS1}
For any arbitrarily small $\ep>0$, one has
\begin{equation}
\label{eq:mod-cont}
\lim_{h \to 0} \sup_{t\in [0, 1]^N, |s-t|\le h}\frac{|\widetilde{Z}^{\vec H}(s)-
\widetilde{Z}^{\vec H}(t)|} {\sum_{j=1}^N |s_j-t_j|^{H_j}\big|\log
\big(\sum_{j=1}^N|s_j-t_j|^{H_j}\big) \big|^{(2 +\ep)/\a}}= 0,
\qquad \hbox{ a.s.}
\end{equation}
\end{theorem}

\begin{proof}\ Let  $\rho$ be the metric on $\R^N$ defined in (\ref{Def:rho})
and let $\{\widetilde D_n, n \ge 1\}$ be the sequence defined in (\ref{Eq:Dn2}).
We claim that for all integers $n \ge 1$, $0 < \ga < \min\{1, \a\}$ and $\eta> 0$,
\begin{equation}\label{Eq:HFSS-1}
\E\bigg(\max_{\tau_n \in \widetilde D_n} \max_{\tau'_{n-1} \in O_{n-1}(\tau_{n})}
\big|\widetilde{Z}^{\vec H}(\tau_n) - \widetilde{Z}^{\vec H}(\tau_{n-1}')\big|^\gamma\bigg)
\le K\, 2^{-\gamma n}\,n^{\frac {(1+\eta)\ga} \a}.
\end{equation}

Before proving (\ref{Eq:HFSS-1}), we note that it implies
\[
\sum_{p=n}^\infty \E\bigg(\max_{\tau_p \in \widetilde{D}_p}
\max_{\tau'_{p-1}\in O_{p-1}(\tau_p)}
\big|\widetilde{Z}^{\vec H} (\tau_p) - \widetilde{Z}^{\vec H}(\tau_{p-1}')\big|^\gamma\bigg)
\le K\, 2^{-\ga n } \,n^{\frac {(1+\eta)\ga} \a}.
\]
Hence $\widetilde{Z}^{\vec H}$ satisfies the conditions of  Corollary \ref{Thm:MC2} with
$\sigma(h) = h \big(\log 1/h\big)^{(1 + \eta)/\a}$. Thus, (\ref{eq:mod-cont}) follows from
(\ref{Eq:UM1d}).

It remains to prove (\ref{Eq:HFSS-1}). Since this is similar to the proof
of Proposition \ref{Prop:HarmF1}, we only provide a sketch of it.
Let $\eta$ and $\beta$ be two positive constants with $\beta$ satisfying
\[
\max_{1 \le j \le N} \bigg\{1 - \frac{2 \a (1 - H_j)} {2 - \a}\bigg\} < \beta < 1.
\]
Define the density function $\varphi: \R^N
\to [0, \infty)$  by $\varphi(x) = \prod_{j=1}^N \varphi_j(x_j)$, where
$x = (x_1, \ldots, x_N)$ and
\[
\varphi_j(x_j) = \left\{\begin{array}{ll}
K_{15}\, |x_j|^{-\beta}\qquad \qquad \qquad \qquad &\hbox{ if }\ |x_j| \le 3,\\
K_{16}\, \big(|x_j| (\log |x|)^{1 + \eta}\big)^{-1}
&\hbox{ if }\ |x_j| > 3.
\end{array}
\right.
\]
In the above, the constants $K_{15}$ and $K_{16}$ are chosen such that
$\int_{\R^N}\varphi(x)\,dx = 1$. By Lemma \ref{Lem:Lepage},
$\widetilde{Z}^{\vec H}$ has the same finite dimensional
distributions as
\[
Y(t) = C_\a\, \Re\bigg(\sum_{j=1}^\infty \Gamma_j^{-1/\alpha}
\varphi(\xi_j)^{-1/\alpha}\, h(t, \xi_j)\, g_j\bigg),
\]
where the function $h(t, x)$ is now defined by
$$
h(t, x) = \prod_{j=1}^N  \frac{e^{i t_j x_j} - 1} {|x_j|^{H_j + \frac 1 \a}}.
$$

Again, by conditional on $\{(\xi_j, \Gamma_j),\, j \ge 1\}$, we have
\begin{equation}\label{Eq:IncVar10}
\E_g \Big[\big(Y(t) - Y(s)\big)^2 \Big] =  C_\a^2\,
\sum_{j=1}^\infty \Gamma_j^{-2/\alpha}
\varphi(\xi_j)^{-2/\alpha}\, \big|h(t, \xi_j) - h(s, \xi_j) \big|^2.
\end{equation}
If follows from (\ref{Eq:IncVar10}) and Lemma \ref{Lem:GaussianMI}
that
\begin{equation}\label{Eq:IncVar11}
\begin{split}
&\E\bigg(\max_{\tau_n \in \widetilde D_n} \max_{\tau'_{n-1} \in O_{n-1}(\tau_{n})}
\big|Y(\tau_n) - Y(\tau_{n-1}')\big|^\gamma\bigg)\\
&\le K \, \, n^{\ga/2} \E_{\Gamma, \xi}\Bigg\{ \bigg(\sum_{j=1}^\infty \Gamma_j^{-2/\alpha}
\varphi(\xi_j)^{-2/\alpha}\, \max_{\tau_n \in \widetilde D_n} \max_{\tau'_{n-1} \in O_{n-1}(\tau_{n})}
\big|h(\tau_n, \xi_j) - h(\tau'_{n-1}, \xi_j) \big|^2\bigg)^{\ga/2}\Bigg\}.
\end{split}
\end{equation}
It can be verified that for every $s, t$ and $x \in \R^N$
\begin{equation}\label{Eq:IncVar13}
\begin{split}
&\big|h(t, x) - h(s, x) \big|^2 = \bigg|\prod_{k=1}^N \big(e^{i t_k x_k} - 1\big)
- \prod_{k=1}^N \big(e^{i s_k x_k} - 1\big)\bigg|^2 \prod_{k=1}^N  \frac{1}
{|x_k|^{2H_k + \frac 2 \a}}\\
&\le K \,\sum_{\ell = 1}^N\Bigg[ \big\{1 \wedge (t_\ell- s_\ell)^2 x_\ell^2\big\}\,
\bigg(\prod_{k < \ell} \big\{1 \wedge t_k^2 x_k^2\big\}\, \prod_{k > \ell}
\big\{1 \wedge s_k^2 x_k^2\big\}\bigg)\Bigg] \prod_{k=1}^N  \frac{1}
{|x_k|^{2H_k + \frac 2 \a}}.
\end{split}
\end{equation}
By using (\ref{Eq:IncVar13}) and an argument similar to the proof of Proposition
\ref{Prop:HarmF1}, one can derive
\[
\begin{split}
&\E_{\Gamma, \xi}\Bigg\{\bigg(\sum_{j=1}^\infty \Gamma_j^{-2/\alpha}
\varphi(\xi_j)^{-2/\alpha}\, \max_{\tau_n \in \widetilde D_n} \max_{\tau'_{n-1}
\in O_{n-1}(\tau_{n})}
\big|h(\tau_n, \xi_j) - h(\tau'_{n-1}, \xi_j) \big|^2\bigg)^{\ga/2}\Bigg\} \\
& \le K\,  2^{-\ga n}\, n^{-\ga(1 + \eta)(\frac 1 2 - \frac \ga \a)}.
\end{split}
\]
This, together with (\ref{Eq:IncVar11}), proves (\ref{Eq:HFSS-1}).
\end{proof}

\vspace{.2in}

\bibliographystyle{plain}

\begin{thebibliography}{1234}

\bibitem{AT07}
Adler, R. J. and Taylor, J. E. (2007), {\it Random Fields and Geometry}.
Springer, New York.

\bibitem{Anh1}
V. V. Anh, J. M. Angulo and M. D. Ruiz-Medina (1999), Possible
long-range dependence in fractional random fields. {\it J. Statist.
Plann. Inference} {\bf 80}, 95--110.

\bibitem{ARX07a}
Ayache, A., Roueff, F. and Xiao, Y. (2007),  Local and asymptotic
properties of linear fractional stable sheets. {\it C. R. Acad.
Sci. Paris, Ser. A.} {\bf 344}, 389--394.

\bibitem{ARX07b}
Ayache, A., Roueff, F. and Xiao, Y. (2009), Linear fractional stable
sheets: wavelet expansion and sample path properties. {\it Stoch.
Process. Appl.} {\bf 119}, 1168–-1197.

\bibitem{Bernard70}
Bernard, P. (1970), Quelques propri\'et\'es des trajectoires des
fonctions al\'eatoires stables sur $\R^{k}$. {\it Ann. Inst. H.
Poincar\'e Sect. B (N.S.)} {\bf 6}, 131--151.

\bibitem{BiermeLacaux07}
Bierm\'e, H. and Lacaux, C. (2009), H\"older regularity for operator
scaling stable random fields. {\it Stoch. Process. Appl.} to appear.

\bibitem{csaki-csorgo}
Cs\'{a}ki, E. and Cs\"org\H o, M. (1992), Inequalities for increments of
stochastic processes and moduli of continuity.
\emph{Ann. Probab.} {\bf 20}, 1031--1052.




\bibitem{Garsia71}
Garsia, A. M. (1972), Continuity properties of Gaussian processes
with multidimensional time parameter. {\it Proc. 6th Berkeley
Symp. Math. Statist. Probability} \textbf{II}, pp. 369--374,
University of California Press, Berkeley.





\bibitem{KhD02}
Khoshnevisan, D. (2002),  {\it Multiparameter Processes: An
Introduction to Random Fields}. Springer, New York.

\bibitem{KoMae91a}
K\^ono, N. and Maejima, M. (1991a), Self-similar stable processes
with stationary increments. In {\it Stable processes and related
topics (Ithaca, NY, 1990)}, 275--295, Progr. Probab., 25,
Birkh\"auser Boston, Boston, MA.

\bibitem{KoMae91b}
K\^ono, N. and Maejima, M. (1991b), H\"older continuity of sample
paths of some self-similar stable processes. {\it Tokyo J. Math.}
{\bf 14}, 93--100.

\bibitem{KwapienR04}
Kwapie\'n, S. and Rosi\'nski, J. (2004), Sample H\"older continuity of
stochastic processes and majorizing measures. In: {\it Seminar on
Stochastic Analysis, Random Fields and Applications IV}, pp. 155--163,
Progr. Probab., {\bf 58}, Birkh\"auser, Basel.


\bibitem{MPisier84}
Marcus, M. B. and Pisier, G. (1984), Characterizations of almost
surely continuous $p$-stable random Fourier series and strongly
stationary processes. {\it Acta Math.} {\bf 152}, 245--301.


\bibitem{MarcusRosen06}
Marcus, M. B. and Rosen, J. (2006), {\it Markov Processes, Gaussian Processes,
and Local Times}. Cambridge University Press, Cambridge.

\bibitem{MRosinski05}
Marcus, M. B. and Rosinski, J. (2005), Continuity and boundedness
of infinitely divisible processes: a Poisson point process approach.
{\it J. Theoret. Probab. } {\bf 18}, 109--160.



\bibitem{Nolan89}
Nolan, J. (1989), Local nondeterminism and local times for stable processes.
{\it Probab. Th. Rel. Fields} {\bf 82}, 387--410.



\bibitem{Pitman68}
Pitman, E. J. G. (1968), On the behavior of the characteristic function of
a probability sidtribution in the neighbourhood of the origin.
{\it J. Australian Math. Soc. Series A} {\bf 8}, 422--443.


\bibitem{Rosinski89}
Rosi\'nski, J. (1989), On path properties of certain infinitely divisible processes.
{\it Stochastic Proc. Appl.} {\bf 33}, 73--87.

\bibitem{Rosinski90}
Rosi\'nski, J. (1990), On series representations of infinitely divisible random
vectors. {\it Ann. Probab.} {\bf 18}, 405--430.


\bibitem{Samoro04}
Samorodnitsky, G. (2004), Extreme value theory, ergodic
theory and the boundary between short memory and long
memory for stationary stable processes. {\it Ann. Probab.}
{\bf 32}, 1438--1468.

\bibitem{ST94}
Samorodnitsky, G. and  Taqqu, M. S. (1994),  {\it Stable
non-Gaussian Random Processes:  Stochastic models with infinite
variance.} Chapman \& Hall, New York.



\bibitem{Taka89}
Takashima, K. (1989), Sample path properties of ergodic
self-similar processes. {\it Osaka Math. J.} {\bf 26}, 159--189.



\bibitem{Tal06}
Talagrand, M. (2006), {\it Generic Chaining}. Springer-Verlag, New York.

\bibitem{Xiao06}
Xiao, Y. (2006), Properties of local nondeterminism of Gaussian and
stable random fields and their applications.  {\it Ann. Fac. Sci.
Toulouse Math.} {\bf XV}, 157--193.

\bibitem{Xiao08}
Xiao, Y. (2008), Properties of strong local nondeterminism and local times of
stable random fields. {\it Seminar on Stochastic Analysis,
Random Fields and Applications VI}, to appear.

\end{thebibliography}
\begin{small}

\end{small}

\vspace{.2in}

\begin{quote}
\begin{small}

\noindent \textsc{Yimin Xiao}.\
        Department of Statistics and Probability, A-413 Wells
        Hall, Michigan State University,
        East Lansing, MI 48824, U.S.A.\\
        E-mail: \texttt{xiao@stt.msu.edu}\\
        URL: \texttt{http://www.stt.msu.edu/\~{}xiaoyimi}

\end{small}
\end{quote}
\end{document}